\documentclass[12pt]{amsart}
\usepackage[mathscr]{eucal}
\usepackage{amssymb}
\usepackage{amsfonts}

\theoremstyle{plain} 
\newtheorem{thm}{Theorem}[section]
\newtheorem{cor}[thm]{Corollary}
\newtheorem{lem}[thm]{Lemma}
\newtheorem{prop}[thm]{Proposition}

\theoremstyle{definition}
\newtheorem{defn}[thm]{Definition}
\newtheorem{ex}[thm]{Example}

\theoremstyle{remark}
\newtheorem{rem}[thm]{Remark}

\def\N{{\mathbb N}}
\def\ZZ{{\mathbb Z}}

\def\P{P_\bullet}

\def\Pc{P^\bullet}
\def\X{X^\bullet}
\def\Y{Y^\bullet}
\def\Z{Z^\bullet}

\def\gr{\operatorname{*mod}}
\def\Gr{\operatorname{*Mod}}
\def\qf{\operatorname{qf}}
\def\fp{\operatorname{*fp}}
\def\cC{{\mathcal C}}
\def\cD{{\mathcal D}}

\def\too{\longrightarrow}

\def\lin{\operatorname{lin}}
\def\Hom{\operatorname{Hom}}
\def\uHom{\underline{\Hom}}
\def\Ext{\operatorname{Ext}}
\def\uExt{\underline{\Ext}}
\def\ld{\operatorname{ld}}
\def\reg{\operatorname{reg}}
\def\Sym{\operatorname{Sym}}
\def\Im{\operatorname{Im}}
\def\soc{\operatorname{soc}}
\def\Ker{\operatorname{Ker}}

\def\cF{{\mathcal F}}
\def\cS{{\mathcal S}}

\def\cH{\mathcal H}
\def\m{{\frak m}}

\def\pd{\operatorname{proj.dim}}

\def\<{{\langle}}
\def\>{{\rangle}}
\def\op{{\sf op}}

\title[{Linearity Defect and Regularity over a Koszul algebra}]
{Linearity Defect and Regularity \\over a Koszul algebra}

\author{Kohji Yanagawa}
\address{Department of Mathematics, Faculty of Engineering Science, 
Kansai University,  Suita 564-8680, Japan}
\email{yanagawa@ipcku.kansai-u.ac.jp}

\begin{document}

\maketitle
\thispagestyle{empty}

\begin{abstract}
Let $A = \bigoplus_{i \in \N}A_i$ be a Koszul algebra over a field $K = A_0$, 
and $\gr A$ the category of finitely generated graded left $A$-modules.   
The {\it linearity defect} $\ld_A(M)$ of $M \in \gr A$ is an invariant 
defined by Herzog and Iyengar. 
An exterior algebra $E$ is a Koszul algebra which is the Koszul dual 
of a polynomial ring. Eisenbud et al. showed that  
$\ld_E(M) < \infty$ for all $M \in \gr E$. Improving their result, 
we show that the  Koszul dual $A^!$ of a Koszul {\it commutative} algebra  
$A$ satisfies the following. 
\begin{itemize}
\item Let $M \in \gr A^!$.  If $\{ \,  \dim_K M_i \mid i \in  \ZZ \, \}$ is bounded,  
then $\ld_{A^!}(M) < \infty$. 
\item If $A$ is complete intersection, then $\reg_{A^!}(M) < \infty$ and 
$\ld_{A^!} (M) < \infty$ for all $M \in \gr A^!$. 
\item If $E=\bigwedge \<y_1, \ldots, y_n\>$ is an exterior algebra, 
then $\ld_E(M) \leq c^{n!} 2^{(n-1)!}$ for $M \in \gr E$    
with  $c := \max \{ \, \dim_K M_i \mid i \in \ZZ \, \}$.
\end{itemize}
\end{abstract}

\section{Introduction}
Let $A =\bigoplus_{i \in \N} A_i$ be a (not necessarily commutative) graded algebra over a 
field $K:= A_0$ with $\dim_K A_i < \infty$ for all $i \in \N$, 
and $\gr A$ the category of finitely generated 
graded left $A$-modules. Throughout this paper, we assume that $A$ is {\it Koszul}, that is, 
$K = A/\bigoplus_{i \geq 1} A_i$ has a graded free resolution of the form 
$$\cdots \too A(-i)^{\beta_i(K)} \too \cdots \too 
A(-2)^{\beta_2(K)} 
\too A(-1)^{\beta_1(K)} \too A \too K \too 0.$$ 
Koszul duality is a certain derived equivalence between 
$A$ and its Koszul dual algebra $A^!:= \Ext^\bullet_A(K,K)$. 

For $M \in \gr A$, we have its minimal graded free resolution 
$\cdots \to P_1 \to P_0 \to M \to 0$, and natural numbers    
$\beta_{i, \, j}(M)$ such that $P_i \cong \bigoplus_{j \in \ZZ} A(-j)^{\beta_{i, \, j}(M)}$. 
We call 
$$\reg_A (M) := \sup \{ \, j -i \mid i \in \N, j \in  \ZZ \ 
\text{with $\beta_{i, \, j}(M) \ne 0$} \, \}$$
the {\it regularity} of $M$. 
If $A$ is not left noetherian, then there is some $M \in \gr A$ such that 
$\sum_{j \in \ZZ} \beta_{1,\,j}(M) = \infty$. In this case, $\reg_A(M) = \infty$. 

When $A$ is a polynomial ring, $\reg_A(M)$ is called the {\it Castelnuovo-Mumford regularity} of $M$,  
and has been deeply studied  from both geometric and computational interest. 
Even for a general Koszul algebra $A$, $\reg_A(M)$ is  still 
an interesting invariant closely related to Koszul duality 
(see Theorem~\ref{reg vs Koszul} below). 

Let $\P$ be a minimal graded free resolution of $M \in \gr A$.
The {\em linear part} $\lin (\P)$ of $\P$ is the chain complex such that
$\lin (\P)_i = P_i$ for all $i$ and its differential maps are given 
by erasing all the entries of degree $\ge 2$ from the matrices representing the 
differentials of $\P$.  
According to Herzog-Iyengar \cite{HI}, we call 
$$\ld_A(M) := \sup \{ \, i \mid H_i (\lin(\P)) \ne 0 \, \}$$
the {\it linearity defect} of $M$. 
This invariant is related to the regularity via 
Koszul duality (see Theorem~\ref{ld vs reg} below).  

In \S4, we mainly treat a Koszul {\it commutative} algebra $A$ or its dual $A^!$. 
Even in this case, it can occur  that 
$\ld_A(M) = \infty$ for some $M \in \gr A$ (c.f. \cite{HI}), while 
Avramov-Eisenbud \cite{AE} showed that $\reg_A(M) < \infty$ 
for all $M \in \gr A$.  On the other hand, Herzog-Iyengar \cite{HI} 
proved that if $A$ is complete 
intersection or Golod then $\ld_A(M) < \infty$ for all $M \in \gr A$.  
Initiated by these results, we will show the following. 

\bigskip

\noindent{\bf Theorem~A.} {\it  
Let $A$ be  a Koszul commutative algebra (more generally, a Koszul algebra with 
$\reg_A(M) < \infty$ for all $M \in \gr A$).  
Then we have;

{\rm (1)} Let $N \in \gr A^!$. If $\reg_{A^!}(N) < \infty$ (e.g. $\dim_K N < \infty$), 
              then $\ld_{A^!}(N) < \infty$. 

{\rm (2)} The following conditions are equivalent. 
\begin{itemize}
\item[(a)] $\ld_A (M) < \infty$ for all $M \in \gr A$.
\item[(a')] $\ld_A (M) < \infty$ for all $M \in \gr A$ with $\dim_K M < \infty$.
\item[(b)] If $N \in \gr A^!$ has a finite presentation, 
then $\reg_{A^!} (N) < \infty$. 
\end{itemize}
}

\bigskip

In Theorem~A (2), the implications $(a) \Rightarrow (a') \Leftrightarrow (b)$ 
hold for a general Koszul algebra. 


When $A$ is commutative, B\o gvad and  Halperin \cite{BH} showed that  
$A^!$ is noetherian if and only if $A$ is complete intersection. 
Moreover, by Backelin and Roos \cite[Corollary~2]{BR}, 
if $A$ is a Koszul complete intersection then $\reg_{A^!}(N) < \infty$ 
for all $N \in \gr A^!$.  (Since $A^!$ admits a {\it balanced dualizing complex}, we can explain 
this also by \cite{Jor04}.)
So, in this case, we have $\ld_A(M) < \infty$ for all $M \in \gr A$  
by Theorem~A (2). This is a part of the above result of Herzog and Iyengar. 
Their proof takes slightly different approach, but is also based on 
a similar result in \cite{BR}.   

Let $\fp A^!$ be the full subcategory of $\gr A^!$ consisting of finitely presented modules.  

\bigskip

\noindent{\bf Theorem~B.} {\it  
If $A$ is a Koszul algebra such that $\ld_A(M) < \infty$ for all 
$M \in \gr A$, then $A^!$ is left coherent (in the graded context), 
and $\fp A^!$ is an abelian category. 
If further $A$ is commutative, then Koszul duality gives  
$$\cD^b(\gr A) \cong \cD^b(\fp A^!)^\op.$$   
}

\medskip

\noindent{\bf Corollary~C.} 
{\it Let $A$ be a Koszul commutative algebra. 
If $A$ is Golod, then we have 
$\cD^b(\gr A) \cong \cD^b(\fp A^!)^\op.$ 
If $A$ is a complete intersection, then we have 
$\cD^b(\gr A) \cong \cD^b(\gr A^!)^\op.$}

\bigskip

Let $E:= \bigwedge \<y_1, \ldots, y_n \>$  be an exterior algebra. 
Eisenbud et al. \cite{EFS} 
showed that $\ld_E(N) < \infty$ for all $N \in \gr E$ 
(now this is a special case of Theorem~A, since $E$ is the Koszul dual of a polynomial 
ring $S:= K[x_1, \ldots, x_n]$).  
If $n \geq 2$, then $\sup\{ \, \ld_E(N) \mid N \in \gr E\, \}= \infty$. 
On the other hand, we will see that 
\begin{equation} \label{bound}
\ld_E(N) \leq c^{n!} 2^{(n-1)!} \quad (c := \max \{ \, \dim_K N_i \mid i \in \ZZ \, \})
\end{equation}
for $N \in \gr E$.  

To prove this, we use (a special case of ) a result of 
Brodmann and Lashgari (\cite[Theorem~2.6]{BL})  
stating that if a submodule $M \subset S^{\oplus c}$ 
is generated by elements of degree 1 then  $\reg_S(M) < c^{n!} 2^{(n-1)!}$. 
But a computer experiment suggests that the bound \eqref{bound} could be very far from sharp. 
For example, if $I \subset E$ is a {\it monomial} ideal then we have $\ld_E(E/I) \leq \max \{  n-2, 1  \}$ 
(\cite{OY}). This does not hold for general graded ideals. 
We have a graded ideal $I \subset E$ with $n = 6$ and $\ld_E(E/I) =9$. 
It is not hard to find similar examples, but these are still much lower 
than the value given in \eqref{bound}. 
 
\section{Koszul Algebras and Koszul Duality}
Let $A = \bigoplus_{i \in \N} A_i$ be a graded algebra over a field 
$K:= A_0$ with $\dim_K A_i < \infty$ for all $i \in \N$, 
$\Gr A$ the category of graded left $A$-modules,  
and $\gr A$ the full subcategory of $\Gr A$ consisting of 
finitely generated modules. We say $M  = 
\bigoplus_{i \in \ZZ} M_i \in \Gr A$ is {\it quasi-finite}, if $\dim_K M_i < \infty$ 
for all $i$ and $M_i = 0$ for $i \ll 0$. If $M \in \gr A$, 
then it is clearly quasi-finite. 
We denote the full subcategory of $\Gr A$ consisting of 
quasi-finite modules by $\qf A$. 
Clearly, $\qf A$ is an abelian category with enough projectives. 
For $M \in \Gr A$ and $j \in \ZZ$, 
$M(j)$ denotes the shifted module of $M$ with $M(j)_i = M_{i+j}$.  
For $M,N \in \Gr A$, set $\uHom_A(M,N):= \bigoplus_{i \in \ZZ} 
\Hom_{\Gr A}(M,N(i))$ to be a graded $K$-vector space with 
$\uHom_A(M,N)_i = \Hom_{\Gr A}(M,N(i))$. Similarly, we also define 
$\uExt^i_A(M,N)$. 

Let $\cC(\qf A)$ be the homotopy category of cochain 
complexes in  $\qf A$, and $\cC^-(\qf A)$ 
its full subcategory consisting of complexes which are bounded above 
(i.e., $\X \in \cC(\qf A)$ with $X^i = 0$ for $i \gg 0$). 
We say $\Pc \in \cC^-(\qf A)$ is a free resolution 
of $\X \in \cC^-(\qf A)$, if  each $P^i$ is a free module and  
there is a quasi-isomorphism $\Pc \to \X$. 
We say a free resolution $\Pc$ is 
{\it minimal}, if $\partial(P^i) \subset \m P^{i+1}$ for all $i$. 
Here $\partial$ denotes the differential map, and 
$\m :=\bigoplus_{i > 0} A_i$ is the graded maximal ideal. 
Any $\X \in \cC^-(\qf A)$ has a minimal free resolution, 
which is unique up to isomorphism.  

Regard $K = A/\m$ as a graded left $A$-module, and set 
$$\beta^i_j(\X) := \dim_K \uExt_A^{-i}(\X, K)_{-j} 
\quad \text{and} \quad 
\beta^i(\X) := \sum_{j \in \ZZ}\beta^i_j(\X)$$ 
for $\X \in \cC^-(\qf A)$ and $i,j \in \ZZ$. In this situation,    
if $\Pc \in \cC^-(\qf A)$ is a minimal free resolution of $\X$, then 
we have $P^i \cong \bigoplus_{j \in \ZZ}A(-j)^{\beta^i_j(\X)}$ 
for each $i \in \ZZ$. It is easy to see that 
$\beta^i_j(\X) < \infty$ for each $i,j$. 

Following the usual convention, we often describe 
(the invariants of) a free resolution of a module $M \in \qf A$  
in the homological manner. So we have 
$\beta_{i,j}(M) = \beta^{-i}_j(M)$, and  
a minimal free resolution of 
$M$ is of the form 
$$\P: \ \cdots \too \bigoplus_{j \in \ZZ}A(-j)^{\beta_{1,j}(M)} 
\too \bigoplus_{j \in \ZZ}A(-j)^{\beta_{0,j}(M)} \too M \too 0.$$

We say $A$ is {\it Koszul}, if $\beta_{i, \, j}(K) \ne 0$ 
implies $i=j$, in other words,  
$K$ has a graded free resolution of the form 
$$\cdots \too A(-i)^{\beta_i(K)} \too \cdots \too 
A(-2)^{\beta_2(K)} 
\too A(-1)^{\beta_1(K)} \too A \too K \too 0.$$ 
Even if we regard $K$ as a right $A$-module, we get the  
equivalent definition.

The polynomial ring $K[x_1, \ldots, x_n]$ and 
the exterior algebra $\bigwedge \< y_1, \ldots, y_n \>$ 
are primary examples of Koszul algebras. Of course, there are 
many other important examples.  
In the noncommutative case, many of them are not 
left (or right) noetherian. 
In the rest of the paper, we assume that $A$ is Koszul.

\medskip

{\it Koszul duality} is a derived equivalence between a Koszul algebra 
$A$ and its dual $A^!$. A standard reference of this subject is 
Beilinson et al. \cite{BGS}. But, in the present paper, we follow the convention of 
Mori~\cite{Mori}. 

Recall that Yoneda product makes  
$A^!:=  \bigoplus_{i \in \N} \Ext_A^i(K,K)$ a graded $K$-algebra. 
(In the convention of \cite{BGS}, $A^!$ denotes the opposite algebra 
of our $A^!$. So the reader should be careful.)  
If $A$ is Koszul, then so is $A^!$ and we have $(A^!)^! \cong A$.  
The Koszul dual of the polynomial ring $S:= K[x_1, \ldots, x_n]$ 
is the exterior algebra $E:= \bigwedge \<y_1, \ldots, y_n \>$. 
In this case, since $S$ is regular and noetherian, 
Koszul duality is very simple. It gives an equivalence 
$\cD^b(\gr S) \cong \cD^b(\gr E)$ of the bounded derived categories. 
This equivalence is sometimes called 
{\it Bernstein-Gel'fand-Gel'fand correspondence} 
({\it BGG correspondence} for short). In the general case, 
the description of Koszul duality is slightly technical. 
For example, if $A$ is not left noetherian, 
then $\gr A$ is not an abelian category.  
So we have to treat $\qf A$.  
  
Let $\cC^\uparrow(\qf A)$ be the full subcategory of $\cC(\qf A)$ 
(and $\cC^-(\qf A)$) consisting of complexes $\X$ satisfying 
$$X^i_j = 0 \quad \text{for $i \gg 0$ or $i+j \ll 0$.}$$ 
And let $\cD^\uparrow(\qf A)$ be 
the localization of $\cC^\uparrow(\qf A)$ at quasi-isomorphisms. 
By the usual argument, we see that $\cD^\uparrow(\qf A)$ is equivalent to 
the full subcategory of the derived category $\cD(\qf A)$ 
(and $\cD^-(\qf A)$)
consisting of the complex $\X$ such that 
$$H^i(\X)_j = 0 \quad \text{for $i \gg 0$ or $i+j \ll 0$.}$$     
We also see that $\cD^\uparrow(\qf A)$ is a 
triangulated subcategory of $\cD(\qf A)$. 

We write $V^*$ for the dual space of a $K$-vector space $V$.
Note that if $M \in \Gr A$ then 
$M^* := \bigoplus_{i \in \ZZ} (M_{-i})^*$ is a graded 
right $A$-module. And we fix a basis $\{x_\lambda\}$ 
of $A_1$ and its dual basis $\{y_\lambda\}$ of $(A_1)^* \, 
(= (A^!)_1)$.  Let $(\X, \partial) \in \cC^\uparrow(\qf A)$. 
In this notation, 
we define the contravariant functor 
$F_A: \cC^\uparrow(\qf A) \to \cC^\uparrow(\qf A^!)$ as follows. 
$$F_A(\X)^p_q = \bigoplus A^!_{q+j} \otimes_K 
(X_{-j}^{j-p})^*$$ 
with the differential $d= d'+d''$ given by  
$$d' : A^!_{q+j} \otimes_K (X_{-j}^{j-p})^* \ni a \otimes m
\longmapsto (-1)^p  \sum a y_\lambda  \otimes m x_\lambda \in 
A^!_{q+j+1} \otimes_K (X_{-j-1}^{j-p})^*$$
and 
$$d'' : A^!_{q+j} \otimes_K (X_{-j}^{j-p})^* \ni a \otimes m
\longmapsto a \otimes \partial^*(m)  \in 
A^!_{q+j} \otimes_K (X_{-j}^{j-p-1})^*.$$
The contravariant functor $F_{A^!}: 
\cC^\uparrow(\qf A^!) \to \cC^\uparrow(\qf A)$ is given by a similar way. 
(More precisely, the construction is different, but the result is similar. See the remark below.) 
They induce the contravariant functors 
$\cF_A: \cD^\uparrow(\qf A) \to \cD^\uparrow(\qf A^!)$ 
and $\cF_{A^!}: \cD^\uparrow(\qf A^!) \to \cD^\uparrow(\qf A)$.  

\begin{rem}
In \cite{Mori}, two Koszul duality functors 
are defined individually. 
The functor denoted by $\bar{E}_A$ is the same as our $\cF_A$.  
The other one which is denoted by $\tilde{E}_A$ 
is defined using the operations $\uHom_K(A^!, -)$ and $\uHom_K(-,K)$. 
But, in our case,  it coincides with $F_A$ except the convention of the sign $\pm1$. 
So we do not give the precise definition of $\tilde{E}_A$ here.        
\end{rem}

\begin{thm}[Koszul duality. c.f. \cite{BGS,Mori}]\label{Koszul duality} 
The contravariant functors $\cF_A$ and $\cF_{A^!}$ give an equivalence 
$$\cD^\uparrow(\qf A) \cong \cD^\uparrow(\qf A^!)^\op.$$
\end{thm}

The next result easily follows from Theorem~\ref{Koszul duality} 
and the fact that $\cF_A(K) = A^!$. 

\begin{lem}[{cf. \cite[Lemma~2.8]{Mori}}]
\label{Betti vs Koszul}
For $\X  \in \cD^\uparrow(\qf A)$, we have 
$$\beta^i_j(\X) = \dim H^{-i-j}(\cF_A(\X))_j.$$ 
\end{lem}

\section{Regularity and Linearity Defect}
Throughout this section, 
 $A = \bigoplus_{i \in \N} A_i$ is a Koszul algebra. 

\begin{defn} 
For $\X \in \cD^\uparrow(\qf A)$, we call 
$$\reg_A (\X) := \sup \{ \, i+j \mid i, j \in  \ZZ \ 
\text{with $\beta^i_j(\X) \ne 0$} \, \}$$
the {\it regularity} of $\X$. 
For convenience, we set the regularity of the 0 module 
to be $- \infty$. 
\end{defn}

If $M \in \qf A$ is {\it not} finitely generated, 
then $\beta_{0, \, j}(M) \ne 0$ for arbitrary large 
$j$ and $\reg_A(M) = \infty$. 

If $A$ is a polynomial ring $K[x_1, \ldots, x_n]$ (more generally, 
$A$ is {\it AS regular}), 
then $\reg_A (\X)$ of $\X \in \cD^b(\gr A)$ can be defined in terms of 
the local cohomology modules $H_\m^i(\X)$, see \cite{EG, Jor04, Y7}.  
If $A$ is commutative, it is known that $\reg_A (M) < \infty$ 
for all $M \in \gr A$ (see Theorem~\ref{reg_S} below). 
But this is not true in the non-commutative 
case. In fact, if $A$ is not left noetherian, then $A$ has a graded 
left ideal $I$ which is  not finitely generated, that is,   $\beta_1(A/I) = \beta_0(I) = \infty$. 
In particular, if $A$ is not left noetherian, then $\reg_A (M) = \infty$ for some $M \in \gr A$. 
The author does not know any example $M \in \gr A$ such that $\beta_i(M) < \infty$ 
for all $i$ but $\reg_A(M) = \infty$.

\begin{lem}\label{reg lemma}\begin{itemize}
\item[(1)] For $M \in \qf A$, we have  
$$\reg_A (M) < \infty \ \Rightarrow \ \beta_i(M) < \infty \ 
\text{for all $i$} \ \Rightarrow 
\text{$M$ has a finite presentation.}$$ 
\item[(2)] If $\X \to \Y \to \Z \to \X[1]$ is a triangle in 
$\cD^\uparrow(\qf A)$, then we have 
$$\reg_A(\Y) \leq \max \{ \, \reg_A (\X), \, \reg_A (\Z) \, \}.$$ 
If $\reg_A(\X) \ne \reg_A(\Z)+1$, then equality holds. 
\item[(3)] If $M \in \gr A$ has finite length, then 
$\reg_A (M) \leq \max \{ \, i \mid M_i \ne 0 \}$. 
\item[(4)] For $\X \in \cD^\uparrow(\qf A)$, we have 
$$\reg_A(\X) \leq \sup \{ \, \reg_A(H^i(\X))+i \mid i \in \ZZ  \, \}.$$
\end{itemize}
\end{lem}

\begin{proof}
(1) is clear. Let us prove (2). 
Since the triangle yields the long exact sequence 
$\cdots \to \uExt^i_A(\Z,K) \to  \uExt^i_A(\Y, K) \to \uExt^i_A(\X, K)  
\to  \uExt^{i+1}_A(\Z, K) \to \cdots$, we have 
the assertions.  

We can prove (3) by induction on $\dim_K M$. 
More precisely, if we set  $d:=\max \{ \, i \mid M_i \ne 0 \}$, 
we have a short exact sequence $0 \to K(-d) \to M \to M' \to 0$. 
Now use the induction hypothesis and (2) of this lemma. 

In \cite[Lemma~2.10]{Y7}, (4) is proved using the spectral sequence 
$$E_2^{p,q} = \uExt_A^p(H^{-q}(\X), K) \Longrightarrow \uExt^{p+q}_A(\X, K)$$
under the additional assumption that $A$ is regular, left noetherian, and 
$\X$ is bounded. But these assumptions are clearly irrelevant.  
\end{proof}

The next result directly follows from Lemma~\ref{Betti vs Koszul}. 

\begin{thm}[Eisenbud et al \cite{EFS}, Mori \cite{Mori}]
\label{reg vs Koszul}
For $\X \in D^\uparrow(\qf A)$, we have 
$$\reg_A(\X) = - \inf\{\, i \mid H^i(\cF_A(\X)) \ne 0 \, \}.$$  
\end{thm}  

We say a complex $\X \in \cD^\uparrow(\qf A)$ is {\it strongly bounded}, if 
$\X$ is bounded (i.e., $H^i(\X)=0$ for $i \gg 0$ or $i \ll 0$) 
and $\reg_A(\X) < \infty$. 
Let $\cD^{sb}(\qf A)$ be the full subcategory of 
$\cD^\uparrow(\qf A)$ consisting of strongly bounded complexes.   

\begin{prop}
$\cD^{sb}(\qf A)$ is a triangulated subcategory of $\cD(\qf A)$. 
\end{prop}

\begin{proof}
Easily follows from Lemma~\ref{reg lemma} (2).   
\end{proof}

\begin{prop}\label{equiv sb} 
The (restriction of) functors $\cF_A$ and $\cF_{A^!}$ give an equivalence 
$$\cD^{sb}(\qf A) \cong \cD^{sb}(\qf A^!)^\op.$$
\end{prop}

\begin{proof}
By Theorem~\ref{Koszul duality}, it suffices to show that $\cF_A(\X) \in \cD^{sb}(\qf A^!)$ for 
all $\X \in \cD^{sb}(\qf A)$. Since $\reg_A(\X) < \infty$, 
$\cF_A(\X)$ is bounded by Theorem~\ref{reg vs Koszul}. 
Similarly, $\cF_{A^!} (\cF_A(\X))$ is isomorphic to 
$\X$,  which is bounded, we have $\reg_{A^!}(\cF_A(\X)) < \infty$.  
\end{proof}

Let $(\Pc, \partial) \in \cC^\uparrow(\qf A)$ be a complex of 
{\it free} $A$-modules such that $\partial(P^i) \subset \m P^{i+1}$, 
in other words, $\Pc$ is a minimal free resolution of some 
$\X \in \cC^\uparrow(\qf A)$.
According to \cite{EFS}, we define the {\it linear part} 
$\lin(\Pc)$ of $\Pc$ as follows: 
\begin{itemize}
\item[(1)] $\lin(\Pc)$ is a complex with $\lin(\Pc)^i =P^i$.  
\item[(2)] The matrices representing the differentials of 
$\lin(\Pc)$ are given by ``erasing" all the entries of degree $\ge 2$ (i.e., replacing them by 0)  
from the matrices representing the differentials of $\Pc$. 
\end{itemize}
It is easy to check that $\lin(\Pc)$ is actually a complex.  
But, even if $\P$ is a minimal free resolution of $M \in \qf A$, 
$\lin(\P)$ is not acyclic (i.e., $H_i(\lin(\P)) \ne 0$ 
for some $i > 0$) in general. 

\begin{defn}[Herzog-Iyengar \cite{HI}] 
Let $M \in \qf A$ and $\P$ its minimal graded free resolution. We call 
$$\ld_A(M) := \sup \{ \, i \mid H_i (\lin(\P)) \ne 0 \, \}$$
the {\it linearity defect} of $M$. 
\end{defn}

We say $M \in \gr A$ has a {\it linear (free) resolution} if there is 
some $l \in \ZZ$ such that $\beta_{i, \, j}(M) \ne 0$ 
implies that $j-i=l$. In this case, the minimal free resolution $\P$ 
of $M$ coincides with $\lin(\P)$, and $\ld_A(M) = 0$.
For $M \in \qf A$ with 
$\iota:= \inf \{ \, i \mid M_i \ne 0 \, \}$, 
$M$ has a linear resolution, if and only if $\reg_A(M) = \iota$, 
if and only if $\reg_A(M) \leq \iota$. 
As shown in \cite[Theorem~5.4]{Mori}, we have  
$$\reg_A(M) = \inf \{ \, i \, \mid \text{$M_{\geq i} := 
\bigoplus_{j \geq i} M_j$ has a linear resolution}\}.$$
For $i \in \ZZ$ and $M \in \qf A$,  $M_{\<i\>}$ denotes the 
submodule of $M$ generated by the degree $i$ component $M_i$. 
We say $M \in \qf A$ is {\it componentwise linear}, 
if $M_{\<i\>}$ has a linear resolution for all $i \in \ZZ$. 
For example, if $M$ has a linear resolution, then 
it is componentwise linear. 
To see this, it suffices to show that if 
$M = \bigoplus_{i \geq 0} M_i$ with $M_0 \ne 0$ 
has a linear resolution, 
then so does $M_{\<1\>}$. But this follows from 
the short exact sequence 
$0 \to M_{\<1\>} \to M \to M/M_{\<1\>} \to 0$ and  
Lemma~\ref{reg lemma} (2), since $\reg_A(M/M_{\<1\>}) 
= \reg_A(K^{\oplus \dim_K M_0})= 0 = \reg_A (M)$.    
Note that $M$ can be componentwise linear 
even if it is not finitely generated. For example, 
$\bigoplus_{i \in \N} K(-i)$ is componentwise linear. 

\begin{prop}[c.f. \cite{R02, Y7}]\label{R & Y}
For $M \in \qf A$, the following are equivalent. 
\begin{itemize}
\item[(1)] $M$ is componentwise linear. 
\item[(2)] $\ld_A(M) = 0$.
\end{itemize}
\end{prop}

This result has been proved by R\"omer \cite{R02} and 
the author \cite[Proposition~4.1]{Y7} under the assumption 
that $M$ is finitely generated. But this assumption is not 
important, since for each $j$ the submodule of $M$ 
generated by $\{ \, M_i \mid i \leq j \, \}$ 
is finitely generated. In the proof of \cite[Proposition~4.1]{Y7},    
the author carelessly stated that ``if $M \in \gr A$ has a finite length, 
then $\reg_A(M) = \max \{ \, i \mid M_i \ne 0 \, \}$", 
which is clearly false (e.g., the exterior algebra 
$E= \bigwedge \< y_1, \ldots, y_n \>$ satisfies 
$\reg_E (E) = 0$ while $E_n \ne 0$).  
But the correct statement 
(Lemma~\ref{reg lemma} (3)) is enough for the proof.  

The next result follows easily from Proposition~\ref{R & Y}.

\begin{prop}[c.f. \cite{R02, Y7}]\label{ld & comp}
For $M \in \qf A$, we have 
$$\ld_A(M) = \inf \{ \, i \mid \text{$\Omega_i(M)$ is componentwise 
linear} \, \},$$
where $\Omega_i(M)$ is the $i^{\rm th}$ syzygy of $M$. 
\end{prop}

Clearly, we have $\ld_A (M) \leq \pd_A (M)$. The inequality is strict 
quite often. For example, we have 
$\pd_A (M) = \infty$ and $\ld_A (M) < \infty$ for many $M$. 
On the other hand, sometimes $\ld_A (M) = \infty$. 

The next result connects the linearity defect with the regularity via 
Koszul duality.  For a complex $\X$, $\cH(\X)$ denotes the complex such that 
$\cH(\X)^i = H^i(\X)$ for all $i$ and all differentials are 0. 

\begin{thm}[{cf. \cite[Proposition~3.4 and 
Theorem~4.7]{Y7}}]\label{ld vs reg} 
Let $\X \in \cD^\uparrow(\qf A)$, and $\Pc$ a minimal free 
resolution of $\cF_A(\X) \in \cD^\uparrow(\qf A^!)$. 
Then we have
$$\lin(\Pc) = F_A \circ \cH (\X).$$ 
Hence, for $M \in \qf A$, 
$$\ld_A (M) = \sup \{ \, \reg_{A^!} (H^i(F_A(M))) +i
\mid i \in \ZZ \}.$$
\end{thm}

\begin{proof}
The first assertion has been proved in \cite[Proposition~3.4]{Y7} 
under the assumption that $A$ is selfinjective 
(or, has a finite global dimension), 
but the assumption is clearly irrelevant. 
We also remark that the Koszul duality functors used in 
\cite{Y7} are covariant, and the $K$-dual of our $\cF$. 
But the essentially same proof as \cite{Y7} also works here.  
The second assertion follows from Theorem~\ref{reg vs Koszul}, since 
$F_{A^!} \circ \cH \circ F_A(M)$ is the linear part of the minimal free resolution of $M$.  
\end{proof}

\section{Koszul Commutative Algebras and their Dual}
If $A$ is a Koszul commutative algebra and 
$S := \Sym_K A_1$ is the polynomial ring, then we have $A=S/I$ for a 
graded ideal $I$ of $S$. In this situation,  $A$ is {\it Golod} if and only if 
$I$ has a 2-linear resolution as an $S$-module (i.e., $\beta_{i,j}(I) \ne 0$ implies $j=i+2$), 
see \cite[Proposition~5.8]{HI}.  We say $A$ 
{\it comes from a complete intersection by a Golod map} 
(see \cite{BR,HI}, although they do not use this terminology),  
if there is an intermediate graded ring $R$ with 
$S \twoheadrightarrow R \twoheadrightarrow A$ 
satisfying the following conditions:
\begin{itemize}
\item[(1)] $R$ is a complete intersection. 
\item[(2)] Let $J$ be the graded ideal of $R$ such that $A = R/J$. Then $J$ has a 2-linear 
resolution as an $R$-module.  
\end{itemize} 
If this is the case, $R$ is automatically Koszul (since so is $A$).  
Clearly, if $A$ itself is complete intersection or Golod, then it 
comes from a complete intersection by a Golod map.  

\begin{ex}
Set $S = K[s,t,u,v,w]$ and $A =S/(st, uv, sw)$. 
Then $A$ is neither Golod nor complete intersection, 
but comes from a complete intersection by a Golod map 
(as an intermediate ring, take $S/(st, uv)$). 
\end{ex}

The next result plays a key role in this section.  

\begin{thm}[Avramov-Eisenbud \cite{AE}]\label{reg_S}
Let $A$ be a Koszul commutative algebra, and 
$S := \Sym_K A_1$ the polynomial ring. 
Then we have $\reg_A (M) \leq \reg_S (M) < \infty$ 
for all $M \in \gr A$. 
\end{thm}

On the other hand, even if $A$ is Koszul and commutative, 
$\ld_A (M)$ can be infinite for some $M \in \gr A$, as 
pointed out in \cite{HI}.  
In fact, if $\ld_A (M) < \infty$ then the Poincar\'e series 
$P_M(t) = \sum_{i \in \N} \beta_i(M) \cdot t^i$ is rational. 
But there exists a Koszul commutative algebra $A$ 
such that $P_M(t)$ is not rational for some $M \in \gr A$ 
(c.f. \cite{Roos2}). By Theorem~\ref{A^!} (2) below, if  $A$ admits a module $M \in \gr A$ 
with $\ld_A (M) = \infty$, then we can take such an $M$ under the additional assumption that 
$\dim_K M < \infty$. 

But we have the following.

\begin{thm}[Herzog-Iyengar \cite{HI}]\label{HI}
Let $A$ be a Koszul commutative algebra. 
 If $A$ comes from a complete intersection by a Golod map (e.g., $A$ itself is 
complete intersection or Golod), then 
$\ld_A(M) < \infty$ for all $M \in \gr A$.
\end{thm}

Now we are interested in $\reg_{A^!}(N)$ and $\ld_{A^!}(N)$ 
for a Koszul commutative algebra $A$. 
First, we remark the important fact that   
the categories $\gr A^!$ and $\gr (A^!)^\op$ are equivalent 
in this case. In fact, a graded left $A^!$-module 
has a natural graded right $A^!$-module structure, and 
vice versa  (c.f. \cite[\S3]{HI}).  
In particular, $A^!$ is left noetherian if and only if 
it is right noetherian.

For the next result and its proof, we need a few preparations. 
For a  graded ring $B = \bigoplus_{i \in \N} B_i$, let  
$\fp B$ be the full subcategory of $\gr B$ consisting of finitely presented modules.  
We say $B$ is {\it left graded coherent}, if any finitely generated graded left ideal of $B$ 
has a finite presentation. As is well-known, $B$ is left graded coherent if and only if 
$\fp B$ is an abelian subcategory of $\gr B$.

\begin{thm}\label{A^!}
If $A$ is a Koszul commutative algebra, we have the following. 
\begin{itemize}
\item[(1)] Let $N \in \gr A^!$. If $\reg_{A^!} (N) < \infty$, 
then $\ld_{A^!} (N) < \infty$. 
\item[(2)] The following conditions are equivalent.
\begin{itemize}
\item[(a)] $\ld_A (M) < \infty$ for all $M \in \gr A$.
\item[(a')] $\ld_A (M) < \infty$ for all $M \in \gr A$ with $M = \bigoplus_{i=0,1} M_i$.
\item[(b)] $\reg_{A^!} (N) < \infty$ for all $N \in \fp A^!$. 
\end{itemize}
\item[(3)]  Let $N \in \qf A^!$. If there is some $c \in \N$ 
such that $\dim_K N_i \leq c$ for all $i \in \ZZ$, 
then $\ld_{A^!} (N) < \infty$. 
\end{itemize}
\end{thm}

\begin{proof}
(1) The complex $F_{A^!}(N)$ is always bounded above. 
Hence if $\reg_{A^!} (N) < \infty$ then $H^i(F_{A^!}(N)) \ne 0$ 
for only finitely many $i$ by Theorem~\ref{reg vs Koszul}. 
Thus the assertion follows from Theorems~\ref{ld vs reg} and 
\ref{reg_S}.  

(2) The implication $(a) \Rightarrow (a')$ is clear. 

$(a') \Rightarrow (b)$: First 
assume that $N \in \fp A^!$ has a presentation of the form 
$A^!(-1)^{\oplus \beta_1} \to  A^{! \, \oplus \beta_0} \to N \to 0$. 
Then there is $M \in \gr A$ with $M= \bigoplus_{i=0,1} M_i$ 
such that $F_A(M)$ gives this presentation. 
Since $\ld_A (M) < \infty$, we have $\reg_{A^!} (N) < \infty$ 
by Theorem~\ref{ld vs reg}.  

Next take an arbitrary $N \in \fp A^!$.  For a sufficiently large $s$,  
$N_{\geq s} := \bigoplus_{i \geq s} N_i$ has a presentation of the form 
$A^!(-s-1)^{\oplus \beta_1} \to A^!(-s)^{\oplus \beta_0} \to 
N_{\geq s} \to 0$.  (To see this, consider the short exact sequence  
$0 \to N_{\geq s} \to N \to N/N_{\geq s} \to 0$, and use the fact that 
$\reg_{A^!} (N/N_{\geq s}) < s$.) 
We have shown that $\reg_{A^!} (N_{\geq s}) < \infty$. 
So $\reg_{A^!} (N) < \infty$ by the above 
short exact sequence.

$(b) \Rightarrow (a)$: 
First, we show that $A^!$ is left graded coherent in this case. 
Assume the contrary. Then there is a finitely generated graded left ideal $I \subset A^!$ which 
is not finitely presented. Clearly, $A^!/I$ has a finite presentation, but $\beta_2(A^!/I) 
= \beta_1(I)=\infty$, in particular, $\reg_{A^!}(A^!/I) = \infty$. 
This is a contradiction. 

So $\fp A^!$ is an abelian category.  
Each term of $F_A(M)$ is a finite free $A^!$-module, in particular, 
$F_A(M) \in \cC^-(\fp A^!)$. Hence we have $H^i(F_A(M)) \in \fp A^!$ for all $i$.    
By the assumption, $\reg_{A^!} (H^i(F_A(M))) < \infty$. 
On the other hand, $H^i(F_{A}(M)) \ne 0$ for finitely many $i$ 
by Theorems~\ref{reg vs Koszul} and \ref{reg_S}. 
So the assertion follows from Theorem~\ref{ld vs reg}. 

(3) Let $\cS$ be the set of all graded submodules of $A^{\oplus c}$ 
which are generated by elements of degree 1. 
By Brodmann~\cite{B}, there is some $C \in \N$ 
such that $\reg_A (M) \leq \reg_S (M) < C$ for all $M \in \cS$. 
Here $S$ denotes the polynomial ring $\Sym_K A_1$. 
To prove the assertion,  it suffices to show that 
$\reg_A (H^i(\cF_{A^!}(N))) +i  \leq C$ for all $i$. We may assume that $i=0$. 
Note that $H^0(\cF_{A^!}(N))$ is the cohomology of the sequence 
$$A \otimes_K (N_1)^* \stackrel{\partial^{-1}}\too 
A \otimes_K (N_0)^*  \stackrel{\partial^0}\too A \otimes_K (N_{-1})^*.$$
Since $\Im(\partial^0)(-1)$ is a submodule of $A^{\oplus \dim_K N_{-1}}$ 
generated by elements of degree 1 and $\dim_K N_{-1} \leq c$, we have 
$\reg_A(\Im(\partial^0)) <  C$. 
Consider the short exact sequence 
$$0 \too \Ker(\partial^0) \too A \otimes_K (N_0)^* 
\too \Im(\partial^0) \too 0.$$ Since $\reg_A(A \otimes_K (N_0)^*)=0$, we have 
$\reg_A(\Ker(\partial^0)) \leq  C$. 
Similarly, we have $\reg_A(\Im(\partial^{-1})) <  C$. 
By the short exact sequence 
$$0 \too \Im(\partial^{-1}) \too \Ker(\partial^0) 
\too H^0(\cF_{A^!}(N)) \too 0,$$ we are done.  
\end{proof}

\begin{rem}\label{non-commutative}
In Theorem~\ref{A^!} (2), the implications $(a) \Rightarrow (a') \Leftrightarrow (b)$ 
hold for a general Koszul algebra.  

If $A$ is a (not necessarily commutative) Koszul 
algebra satisfying $\reg_A(M) < \infty$ for all $M \in \gr A$, 
then Theorem~\ref{A^!} (1) and (2) hold for $A$. 
\end{rem}

In \cite[Corollary~3]{BR}, Backelin and Roos showed that 
if $A$ is a Koszul commutative algebra which comes from a complete intersection 
by a Golod map then $A^!$ is left graded coherent. 
Moreover, they actually proved that $\reg_{A^!}(N) < \infty$ for all $N \in \fp A^!$ 
(see \cite[Corollary~2]{BR} and \cite[Lemma~5.1]{HI}). 
So we have $\ld_A (M) < \infty$ for all $M \in \gr A$ by Theorem~\ref{A^!}, that is, 
we get a result of Herzog and Iyengar (Theorem~\ref{HI}).  
Their original proof is essentially based on this line too.    

A deep theory on the Hopf algebra structure of $A^!$ plays a key role in \cite{BR}. 
But, when $A$ is a Koszul complete intersection, 
we have another exposition of the fact that $\reg_{A^!}(N) < \infty$ for all $N \in \gr A^!$. 
Since this exposition has its own interest, we will give it here. 
The next lemma might be known to specialists. But the author could not find 
reference. So we give a proof, which is suggested by Professor Izuru Mori. 
For the unexplained terminology appearing 
in the next result and its proof, consult \cite{Mori,Smi,Ye}.

\begin{lem}\label{Mori}
If $A$ is a complete intersection, then $A^!$ is left noetherian 
and admits a balanced dualizing complex. 
\end{lem}

\begin{proof}
Let  $S:= \Sym_K S_1$ be the polynomial ring.  
Then we have a regular sequence $z_1, \ldots, z_m \in S_2$   
such that $A = S/(z_1, \ldots, z_m)$. 
Recall that $E:=S^!$ is the exterior algebra. 
Set $A_{(1)}:= S/(z_1)$. Then  
there is a central regular element $w_1 \in (A_{(1)})^!$ of degree 2  
such that $(A_{(1)})^!/(w_1) \cong E$ by \cite[Theorem 5.12]{Smi}. 
Since $E$ is artinian, then it is noetherian and admits a 
balanced dualizing complex. 
Hence $(A_{(1)})^!$ is noetherian and admits a 
balanced dualizing complex by \cite[Lemma~7.2]{Mori}. 
Similarly, if we set $A_{(2)} := S/(z_1, z_2) =A_{(1)}/(\bar{z}_2)$, 
then there is a central regular element $w_2 \in (A_{(2)})^!$ of degree 2  
such that $(A_{(2)})^!/(w_2) \cong A_{(1)}$. 
Hence $(A_{(2)})^!$ is noetherian and admits a balanced dualizing 
complex again. Repeating this argument, we see that 
$A^!$ is noetherian and has a balanced dualizing complex.  
\end{proof}

\begin{cor}\label{A^! for c.i.}
If $A$ is a Koszul complete intersection, then $\reg_{A^!} (N) < \infty$ and 
$\ld_{A^!} (N) < \infty$ for all $N \in \gr A^!$. 
\end{cor}

\begin{proof}
By Lemma~\ref{Mori} and \cite{Jor04}, we have 
$\reg_{A^!} (N) < \infty$ for all $N \in \gr A^!$. 
Hence $\ld_{A^!} (N) < \infty$ for all $N \in \gr A^!$ by Theorem~\ref{A^!}. 
\end{proof}

Next we will treat a Koszul algebra $A$ 
such that $\ld_A (M) < \infty$ for all $M \in \gr A$. 
In this case, $\reg_{A^!} (N) < \infty$ for all 
$N \in \fp A^!$ by Remark~\ref{non-commutative}.
So we have the following (see the proof of 
the implication (b) $\Rightarrow$ (a) of Theorem~\ref{A^!} (2)).

\begin{prop}\label{coherent}
Let $A$ be a Koszul algebra.  If $\ld_A (M) < \infty$ 
for all $M \in \gr A$, then $A^!$ is left graded coherent.  
\end{prop}

\begin{lem}\label{finite presentation}
Assume that $\reg_{A^!} (N) < \infty$ for 
all $N \in \fp A^!$.  Let $\X \in \cD^b(\qf A^!)$ be a bounded complex. 
Then $\X$ is strongly bounded if and only if $H^i(\X) \in \fp A^!$ for all $i$. 
\end{lem}

\begin{proof}
(Sufficiency): If $H^i(\X) \in \fp A^!$, then 
$\reg_{A^!} (H^i(\X)) < \infty$. Since $\X$ is bounded, we have 
$\reg_{A^!}(\X) < \infty$ by Lemma~\ref{reg lemma} (4).

(Necessity): Assume that $\X$ is strongly bounded (more generally, 
$\beta^i(\X) < \infty$ for all $i$). 
Let $\Pc$ be a minimal free resolution of $\X$. Clearly, $\Pc \in \cC^-(\fp A^!)$. 
By Proposition~\ref{coherent},  $\fp A^!$ is an abelian category.  
Hence each $H^i(\Pc) \, (\cong H^i(\X))$ belongs to $\fp A^!$.  
\end{proof}

If $A$ is commutative, then $A$ is noetherian and 
$\gr A$ is an abelian category. So we can consider the derived 
category $\cD^b(\gr A)$, which is a full subcategory of 
$\cD^\uparrow(\qf A)$.   

\begin{lem}\label{Koszul duality for commutative}
Let $A$ be a  Koszul commutative algebra. Then 
$\cD^b(\gr A) = \cD^{sb}(\qf A)$ and the Koszul duality gives   
$\cD^b(\gr A) \cong \cD^{sb}(\qf A^!)^\op.$    
\end{lem}

\begin{proof}
By Proposition~\ref{equiv sb}, 
it suffices to show the first statement. If $\X \in \cD^b(\gr A)$, then 
$\reg_A(\X) < \infty$ by Lemma~\ref{reg lemma} (4) and Theorem~\ref{reg_S}.  
Hence we have $\X \in \cD^{sb}(\qf A)$. 
Conversely, if $\Y \in \cD^{sb}(\qf A)$, then $\beta^i(\Y) < \infty$ 
for all $i$, and the minimal free resolution of $\Y$ is a complex 
of finite free modules. So we have $\Y \in \cD^b(\gr A)$.  
\end{proof}

\begin{thm}\label{Golod}
Let $A$ be a Koszul commutative algebra such that $\ld_A(M) < \infty$ 
for all $M \in \gr A$ (e.g. $A$ comes from a complete intersection by a Golod map). 
Then Koszul duality gives an equivalence 
$\cD^b(\gr A) \cong \cD^b(\fp A^!)^\op.$
\end{thm}

\begin{proof}
By Proposition~\ref{coherent},  
$\fp A^!$ is an abelian category, and closed under extensions in $\qf A^!$. 
Since a free $A^!$-module of finite rank belongs to $\fp A^!$,  
this category has enough projectives. So we have $\cD^b(\fp A^!) = 
\cD^b_{\fp A^!}(\qf A^!) = \cD^{sb}(\qf A^!)$. Here the first equality 
follows from \cite[Exercise~III.2.2]{GM} and the second one follows 
from Lemma~\ref{finite presentation}.  
Now the assertion is a direct consequence of 
Lemma~\ref{Koszul duality for commutative}. 
\end{proof}

We remark that the next corollary also follows from 
Lemma~\ref{Mori} and \cite[Proposition~4.5]{Mori}.  

\begin{cor}\label{CI}
If $A$ is a Koszul complete intersection, then Koszul duality gives 
$\cD^b(\gr A) \cong \cD^b(\gr A^!)^\op.$ 
\end{cor}

In the rest of the paper, we study the linearity defect over the exterior algebra 
$E:= \bigwedge \<y_1, \ldots, y_n \>$. Eisenbud et al. \cite{EFS} showed that $\ld_E (N) < \infty$ 
for all $N \in \gr E$. Now this is a special case of 
Theorem~\ref{A^!}. Every part of the theorem induces their result. 
But the behavior of $\ld_E(N)$ is still mysterious.

If $n \geq 2$, then we have $\sup \{ \, \ld_E(N) \mid  N \in \gr E \, \} = \infty$. 
In fact, $N:=E/\soc(E)$ satisfies $\ld_E(N) \geq 1$. 
And the $i^{\rm th}$ {\it cosyzygy} $\Omega_{-i}(N)$ of $N$ 
(since $E$ is selfinjective, we can consider cosyzygies) 
satisfies $\ld_E(\Omega_{-i}(N)) > i$. 
But we have an upper bound of $\ld_E(N)$ depending only on 
$\max \{ \, \dim_K N_i \mid i \in \ZZ \, \}$ and $n$. 
Before stating this, we recall a result on  $\reg_S(M)$ for $M \in \gr S$. 

\begin{thm}[Brodmann and Lashgari, {\cite[Theorem~2.6]{BL}}]\label{BL}
Let $S= k[x_1, \ldots, x_n]$ be the polynomial ring. 
Assume that a graded submodule $M \subset S^{\oplus c}$ is generated by 
elements whose degrees are at most $d$. Then we have 
$\reg_S(M) < c^{n!} (2d)^{(n-1)!}.$ 
\end{thm}

When $c=1$ (i.e., when $M$ is an ideal), the above bound is a classical result, 
and there is a well-known example which shows the bound is rather sharp. 
For our study on $\ld_E(N)$, the case when $d=1$ (but $c$ is general) is essential. 
When $c=d=1$, we have $\reg_S(M) =1$ for all $M \in \gr S$.  
So the author believes the bound can be strongly improved when $d=1$. 

\begin{prop}
Let $E = \bigwedge \< y_1, \ldots, y_n \>$ be an exterior algebra, 
and $N \in \gr E$. Set $c := \max \{ \, \dim_K N_i \mid i \in \ZZ \, \}$. 
Then $\ld_E(N) \leq c^{n!} 2^{(n-1)!}$.   
\end{prop}

\begin{proof} 
If $M$ is a graded submodules of $S^{\oplus c}$ 
generated by elements of degree 1, 
then we have $\reg_S(M) < c^{n!} 2^{(n-1)!}$ by Theorem~\ref{BL}.   
Now the assertion follows from the argument similar to the proof 
of  Theorem~\ref{A^!} (3).   
\end{proof}

\section*{Acknowledgments}
The author is grateful to Professors Izuru Mori and 
Ryota Okazaki for useful comments.

\end{document}